\documentclass[11pt]{article}

\usepackage{graphicx}
\usepackage{amsmath}
\usepackage{amsthm}
\usepackage{amssymb}
\usepackage[mathscr]{eucal}
\usepackage[all]{xy}

\title{The Birman-Craggs-Johnson homomorphism and \\ abelian cycles
in the Torelli group}
\author{Tara E. Brendle and
Benson Farb \thanks{The first author is partially supported by
NSF grant DMS-0504208 and was also supported in part
by a VIGRE postdoc under NSF grant
number 9983660 to Cornell University.  The second author is
supported in part by NSF grant DMS-0244542.}}

\theoremstyle{plain}
\newtheorem{theorem}{Theorem}
\newtheorem{proposition}[theorem]{Proposition}

\newtheorem{corollary}[theorem]{Corollary}

\newtheorem{question}[theorem]{Question}

\def\proof{{\bf {\medskip}{\noindent}Proof. }}

\def\endproof{$\diamond$ \bigskip}

\def\title{\em}

\def\bar{\overline}

\newcommand\F{\mbox{\bf F}}

\newcommand\calA{{\cal A}}

\newcommand\Z{\mbox{\bf Z}}
\newcommand\Q{\mbox{\bf Q}}
\newcommand\Out{\mbox{Out}}

\newcommand\lowerbound{{16g^4 + O(g^3)}}

\newcommand\mgo{{\rm Mod_{g,1}}}

\newcommand\Mod{{\rm Mod}}

\newcommand\ztwo{\F_2}

\DeclareMathOperator\A{{\cal A}}
\DeclareMathOperator\T{{\cal I}_g}
\DeclareMathOperator\Tgo{{\cal I}_{g,1}}
\DeclareMathOperator\kg{{\cal K}_g}
\DeclareMathOperator\kgo{{\cal K}_{g,1}}

\DeclareMathOperator\Sp{Sp}
\DeclareMathOperator\IM{IM}
\DeclareMathOperator\image{image}

\begin{document}
\maketitle
\begin{abstract}
In the 1970's,
Birman-Craggs-Johnson \cite{BC,Jo1}
used Rochlin's invariant for homology $3$-spheres
to construct a remarkable surjective homomorphism $\sigma:\Tgo\to B_3$,
where $\Tgo$ is the Torelli group and
$B_3$ is a certain $\ztwo$-vector space of Boolean (square-free)
polynomials.  By pulling back cohomology classes and evaluating them on
abelian cycles, we construct $\lowerbound$ dimensions
worth of nontrivial elements of
$H^2(\Tgo, \F_2)$ which cannot be detected rationally.
These classes in fact restrict to nontrivial classes in
the cohomology of the subgroup $\kgo<\Tgo$ generated by
Dehn twists about separating curves.

We also use the ``Casson-Morita algebra'' and Morita's integral lift of
the Birman-Craggs-Johnson map restricted to $\kgo$ to give the same
lower bound on $H^2(\kgo, \Z)$.
\end{abstract}
\tableofcontents

\newpage
\section{Introduction}

Let $\Sigma_{g,1}$ denote the compact, oriented surface of genus $g$
with one boundary component.
The {\em mapping class group} $\Mod_{g,1}$ is the group of
isotopy classes of orientation-preserving self-homeomorphisms
$\Sigma_{g,1}$, where both the homeomorphisms and the isotopies are
taken to fix the boundary $\partial\Sigma_{g,1}$ pointwise.  While
versions of our results extend to closed surfaces, we consider this case
because it is somewhat simpler technically.

Algebraic intersection number gives a symplectic form on
$H_1(\Sigma_{g,1},\Z)$.  This form is preserved by the natural action
of $\Mod_{g,1}$.  The {\em Torelli group} $\Tgo$ is defined to be the kernel
of this action.  We then have an exact sequence
\begin{equation}
\label{eq:exact}
1\to \Tgo\to \Mod_{g,1}\to \Sp(2g,\Z)\to 1
\end{equation}

We will also consider the {\em bounding twist group}
$\kgo$, which is the subgroup of $\Tgo$ generated
by Dehn twists about those curves which separate
$\Sigma_{g,1}$.  Johnson found a homomorphism $\tau$ and proved that it
gives the following exact sequence (see \cite{Jo3,Jo5,Jo6}):

$$1\to\kgo\to\Tgo\stackrel{\tau}{\to}\wedge^3H\to 1$$
where $H=H_1(\Sigma_{g,1};\Z)$.

While there has been spectacular progress in understanding
$H^\ast(\Mod_{g,1},\Z)$ (see \cite{MW}), very little is known about
$H^\ast(\Tgo,\Z)$, and even less is known about $H^\ast(\kgo,\Z)$.  For
example, although it was recently shown in \cite{BF} that $\kgo$ is
infinitely generated, it is still not known whether or not
$H^1(\kgo,\Z)$ is finitely generated.
Note that, as follows from (\ref{eq:exact}), for any field $K$
the vector space
$H^\ast(\Tgo,K)$ is a module over $\Sp(2g,K)$.

Akita \cite{Ak} has shown that the algebras
$H^\ast(\T,\Q)$ and $H^\ast(\kg,\Q)$ must be infinite dimensional for $g\geq
7$, although these groups vanish in sufficiently high degrees.
His proof is by contradiction, however, and produces no explicit
classes. One goal is, therefore, to construct explicit classes.
Morita \cite{Mo2} has identified a certain
($\Mod_{g,1}$-invariant) secondary characteristic classes which generates
$H^1(\kgo, \Z)^{\Mod_{g,1}} \cong \Z^2$.

In a series of papers Johnson proved the difficult result:

$$H_1(\Tgo,\Z)\approx \wedge^3H\oplus B_2$$

\noindent
where $B_2$ consists of $2$-torsion; see \cite{Jo4}, \cite{Jo2} for a
summary, or \cite{vdB}.  While the $\wedge^3H$ piece comes
from purely algebraic considerations, the $B_2$ piece is ``deeper'' in
the sense that it is purely topological, and
comes from the Rochlin invariant (see below);
indeed the former appears in $H_1$ of the ``Torelli group'' in
the analogous theory for $\Out(F_n)$, while
the latter does not.

\bigskip
\noindent
{\bf Problem: }
\label{p:g}
Determine the subalgebras of $H^\ast(\Tgo,K)$, for $K=\Q$ and
$K=\F_2$, generated by $H^1(\Tgo,\Z)\otimes K$.

\bigskip
\noindent
{\bf Rational cohomology. }
Johnson proved in \cite{Jo4} that the induced map in cohomology
\begin{equation}
\label{eq:tau:star}
\tau^\ast:H^\ast(\wedge^3H,\Q)\to H^\ast(\Tgo,\Q)
\end{equation}
is an isomorphism in degree one.  Hain \cite{Hai}
determined the image (and kernel)
of $\tau^\ast$ in degree two, as did Sakasai \cite{Sa}
(up to a small unknown) in degree three.

A key tool in the proofs of these results is the
fact that $\tau^\ast$ is a map of $\Sp(2g,\Q)$-modules.  Since
$H^\ast(\wedge^3H,\Q)$ decomposes as a direct sum of easily described
irreducible representations, the calculation of $\tau^\ast$ is greatly
simplified.  Note that since $\kgo$ lies in (indeed equals) the kernel
of $\tau$, this method yields no information about $\kgo$.

\paragraph{The Birman-Craggs-Johnson homomorphism.}  Consider
a fixed Heegaard embedding $h:\Sigma_g\hookrightarrow S^3$ of the closed
surface $\Sigma_g$.  An element of $f\in\T$ then
determines an integral homology $3$-sphere $M_f$
by cutting out a handlebody determined
by $h$ and gluing it back in via $f$.
Birman-Craggs \cite{BC} proved that the map $f\mapsto \mu(M_f)$ taking
$f$ to the Rochlin invariant $\mu$ of $M_f$ actually determines a
homomorphism $\T \rightarrow
\ztwo$. The is called the {\em Birman-Craggs homomorphism} corresponding
to $h$.

Building on their work, Johnson constructed in \cite{Jo1} the surjective
{\em Birman-Craggs-Johnson (BCJ) homomorphism} $\sigma: \Tgo \rightarrow
B_3$, where $B_3$ is the degree $\leq 3$ summand of a certain graded
$\F_2$-algebra $B$ of ``Boolean polynomials'' in $2g$ variables (see
Section~\ref{spforms} below).  The map $\sigma$ encodes all of the
Birman-Craggs homomorphisms in the sense that the kernel of $\sigma$ is
the intersection of the kernels of all Birman-Craggs homomorphisms.
As with Johnson's homomorphism $\tau$, the group
$\Mod_{g,1}$ acts naturally on both $\Tgo$ and on $B_3$, and $\sigma$
respects this action.

Johnson showed that, unlike the case of the map $\tau$, the
restriction of $\sigma$ to $\kgo$ is highly nontrivial.  In fact,
he showed that:
$$
\sigma|_{\kgo}: \kgo \rightarrow B_2
$$
is a surjection onto the degree $\leq 2$ summand $B_2$ of $B$.  Further,
he showed that the $B_3-B_2$ piece of $\image(\sigma)$ is
precisely $\image(\tau)=\wedge^3H$ reduced mod $2$.  For this reason we
concentrate here on the ``truly $2$-torsion'' summand $B_2$.
Our main result is the following.

\begin{theorem}
\label{theorem:main1}
Each of the images of
$$\sigma^\ast:H^2(B_2,\F_2)\to H^2(\Tgo,\F_2)$$
$$(\sigma|_{\kgo})^\ast:H^2(B_2,\F_2)\to H^2(\kgo,\F_2)$$
have dimension at least $\lowerbound$.  Further, this lower bound is
correct to third order: the subspaces of cocycles we describe below
have codimension at most $4g^2$ in the images of $\sigma^\ast$
(resp. ($\sigma|_{\kgo})^\ast$).
\end{theorem}

Note that none of these classes can be detected via the Johnson
homomorphism $\tau$.  Another
aspect of this problem which interested us is the failure of
representation theory in this context.  While
$\sigma^\ast:H^\ast(B_2,\F_2)\to H^\ast(\Tgo,\F_2)$ is a homomorphism of
$\Sp(2g,\F_2)$-algebras, the module $H^\ast(B_2,\F_2)$ no longer
decomposes as a direct sum of irreducible representations.  Further,
the seemingly simple
(characteristic two) modular representation theory needed
to aid in computations (as in the rational case) seems to be beyond
what is currently known.  Thus we are forced to use
more involved topological methods for computations.

\medskip
\noindent
{\bf The Casson-Morita algebra. } Much, but not all, of the situation
for the BCJ homomorphism described above can be lifted to $\Z$
coefficients.  Replacing the Rochlin invariant $\mu$ with the Casson
invariant $\lambda$, one can consider
the map from the Torelli group of the closed surface $\T\to \Z$ given by
$f\mapsto \lambda(M_f)$.  While this map is definitely not a
homomorphism, Morita proved in \cite{Mo1} that it is a homomorphism when
restricted to $\kg$.  Varying over Heegaard embeddings
$\Sigma_g\hookrightarrow S^3$ gives different homomorphisms.  Morita put
these together (passing for technical reasons to the case of one
boundary component) to give a homomorphism $\rho:\kgo\to
\A$ for a certain $\Z$-algebra $\A$ which we call the {\em Casson-Morita
algebra}, which is a kind of
universal receptor for Casson invariants.
This lifts the BCJ picture as
follows (see Proposition \ref{triangle} below):  there exists a
homomorphism
$\mu: \A \rightarrow B_2$ such that the following diagram commutes:

$$\xymatrix{
& \A \ar[d]^{\mu} \\
\kgo \ar[ur]^{\rho} \ar[r]^{\sigma} & B_2  }
$$

We can then use the homomorphism
\begin{equation}
\rho^\ast:H^\ast(\A,\Z)\to H^\ast(\kgo,\Z)
\end{equation}
to obtain nontrivial classes in $H^\ast(\kgo,\Z)$.  An easy argument
shows that our computations in the $\F_2$ case imply the following.

\begin{corollary}
The image of $\rho^\ast$ in $H^2(\kgo,\Z)$
contains an $\Sp(2g,\Z)$-submodule of rank $\lowerbound$.
\end{corollary}

We remark that, apart from allowing us to use $\sigma^\ast$ (which is
defined on all of $\Tgo$),  the use of $\F_2$
coefficients greatly simplifies computations.  For these reasons,
together with our desire for an
explicit topological construction of cycles, we decided not to work
with Morita's homomorphism $\rho$ (and with $Z$ coefficients)
directly.

\medskip
\noindent
{\bf Abelian cycles. } A pair of commuting elements $f,g$ in a group
$\Gamma$ determines a homomorphism $i:\Z^2\to\Gamma$.  The
{\em abelian cycle} $\{f,g\}\in H_2(\Gamma,\Z)$ determined by the pair
$f,g$ is the image of the fundamental class under the homomorphism
$H_2(\Z^2,\Z)\to H_2(\Gamma,\Z)$.  While Theorem \ref{theorem:main1} is
stated for convenience in cohomology, we mostly study homology, since
our method for proving that the cohomology classes
in the image of $\sigma^\ast$ are nontrivial is to
evaluate these classes on abelian cycles; see \S\ref{abeliancycles} below.

\paragraph{Acknowledgements.}We thank Nathan Broaddus, Ken Brown,
Allen Hatcher, Shigeyuki Morita and Peter Sin for their valuable
expertise and comments.  We are
grateful to Vijay Ravikumar for his useful insights regarding
``index-matched'' elements and for working out the case of bounding pair
maps while participating in Cornell University's REU program in 2005
under the supervision of the first author.  We also wish to thank the
other REU students, Tova Brown, Tom Church, Peter Maceli, and Aaron
Pixton, for helpful discussions.

\section{The Birman-Craggs-Johnson homomorphism}

We begin with a brief review of the work of Birman-Craggs and of Johnson
relating the Rochlin invariant of 3-manifolds to certain algebraic
structures associated to surfaces and their Torelli groups.

\subsection{Birman-Craggs homomorphisms and $\Sp$-quadratic forms}
\label{spforms}

Let $M$ be an oriented,
integral homology $3$-sphere endowed with a spin structure. Choose $X$, a spin
4-manifold such that $\partial X = M$ and such that this
restriction to the boundary induces the given spin structure on
$M$.  The {\em Rochlin invariant} $\mu(M)\in \F_2$ of $M$ is given by
$$
\mu(M):= \frac{\sigma(X)}{8} \hskip .15in \textrm{mod} \hskip
.05in 2 $$ where $\sigma(X)$ denotes the signature of $X$; it does not
depend on the choice of the 4-manifold $X$ (see, e.g. \cite{GS}).  Note
that $\mu(S^3) = 0$.

We give the definition of the Birman-Craggs homomorphisms as
reformulated by Johnson for the case of a surface with one boundary
component \cite{Jo5}. A {\em Heegaard embedding} of the surface
$\Sigma_{g,1}$ in $S^3$ is an embedding $h:\Sigma_{g,1} \rightarrow S^3$
such that $h(\Sigma_{g,1}) \subset S$, where $S$ is a Heegaard surface
for $S^3$.  Let $h$ be such a Heegaard embedding, and let $f \in \T$.
Now split $S^3$ along $h(\Sigma_{g,1})$ and reglue via the map $f$.
Specifically, we map a point $x$ in one copy of the interior of
$h(\Sigma_{g,1})$ to $hfh^{-1}(x)$ in the other copy.  The resulting
3-manifold $M(h,f)$ is necessarily a homology 3-sphere.  We can
therefore define the {\em Birman-Craggs homomorphism} $\rho_h: \Tgo
\to\Z_2$ associated to the Heegaard embedding $h$ via the formula
$$\rho_h(f):=\mu(M(h,f))$$

\paragraph{Self-linking forms.}
Given a Heegaard embedding $h$ of the oriented surface $\Sigma_{g,1}$ in $S^3$,
there is a corresponding mod 2 self-linking form $\omega_h: H_1(\Sigma_{g,1},\Z)
\rightarrow \ztwo$ defined on irreducible elements $c \in
H_1(\Sigma_{g,1},\Z)$  via $\omega_h(c) = lk(h(c),h(c)^+)$; here
the latter expression denotes the linking number of a representative of the
homology class $h(c)$ with its positive push-off $h(c)^+$ in
$S^3$.  For simplicity, we let $H = H_1(\Sigma_{g,1}, \ztwo)$.

Johnson showed that all such mod 2 self-linking forms arise as functions
$\omega: H \rightarrow \ztwo$ whose associated bilinear form is just the
usual symplectic intersection pairing on $H$, and vice-versa. In other
words, if $\Omega$ denotes the set of all mod 2 self-linking forms, then
we have $$\Omega = \{ \omega: H \rightarrow \ztwo |
\hskip .15in \omega(a + b) = \omega(a) + \omega(b) + a \cdot b \}
$$
where $a \cdot b$ denotes the intersection form.  In particular,
given any $\omega \in \Omega$, we have $\omega = \omega_h$ for some
Heegaard embedding $h$.  Moreover, Johnson showed that the value of the
Birman-Craggs homomorphism $\rho_h(f)$ is completely determined by the
self-linking form $\omega_h$ and by $f \in \Tgo$ \cite{Jo1}.

\paragraph{The BCJ map.}  Johnson is then able to combine all Birman-Craggs homomorphisms into a single, surjective map $\sigma$ from $\Tgo$ into the vector space of functions $\Omega \rightarrow \ztwo$, defined as follows:$$
\sigma(f)(\omega) = \rho_h(f) = \mu(h,f)
$$ where $\omega = \omega_h$, as above.  We call $\sigma$ the {\em
Birman-Craggs-Johnson} homomorphism, or BCJ homomorphism for short. It
is clear that $\textrm{ker} \sigma = \cap_{h} \textrm{ker} \rho_h$,
where the intersection is taken over all Birman-Craggs homomorphisms.

\subsection{Johnson's formula for $\sigma$}
Let $a_1, b_1, \ldots, a_g, b_g$ be a fixed symplectic basis for
$H$.   Now for each $c \in H$, we can define a map
$\bar{c}: \Omega \rightarrow \ztwo$ via $\omega \mapsto \omega(c)$.
%\begin{eqnarray*}
%\bar{c}: \Omega &\rightarrow& \ztwo \\
%\omega &\mapsto& \omega(c)
%\end{eqnarray*}
We will need two basic facts which follow directly from the
definitions:
\begin{equation}
\bar{a+b} = \bar{a} + \bar{b} + a \cdot b \mbox{\ \ for all $a,b \in
H$}
\end{equation}
\begin{equation}
\label{squarefree}
\bar{c}^2 = \bar{c} \mbox{\ \ for all $c \in H$}
\end{equation}

The ``square-free'' condition in (\ref{squarefree})
gives the space of such elements a
particularly nice structure, which we now describe.

\paragraph{Boolean polynomials.}
We define $B(x_1, \ldots, x_n)$, the ring of {\em Boolean} (or {\em
square-free}) polynomials, to be the quotient of the usual polynomial
ring on the variables $x_1, \ldots, x_n$ with coefficients in $\ztwo$ by
the ideal generated by the relations $x_i^2 = x_i$ for all $i = 1,
\ldots, n$.  In other words $$B(x_1, \ldots, x_n) = \ztwo [x_1, \ldots,
x_n] / < x_i^2 = x_i >$$ Then we let $B_r = B_r(x_1, \ldots, x_n)$
denote the subset of $B(x_1, \ldots, x_n)$ consisting of those
polynomials of degree at most $r$.  For our
purposes, we will take $n = 2g$,
and we will attach a certain topological significance to these $2g$
variables.

We form Boolean polynomials in the obvious way out of
elements $\bar{c}$, with $c \in H$.  Specifically, we take as
abstract variables the $2g$ maps $\bar{a}_1, \ldots, \bar{a}_g,
\bar{b}_1, \ldots \bar{b}_g$.  As above, we denote the ring of Boolean
polynomials of degree at most $r$ in these $2g$
variables simply by $B_r$.  A basis for $B_2$ as a vector space over
$\F_2$ consists of the $2g^2+g+1$ elements
$$
\{ 1, \bar{a}_i, \bar{b}_i, \bar{a}_i \bar{b}_j, \bar{a}_i
\bar{a}_j (i \neq j), \bar{b}_i \bar{b}_j (i \neq j) \} .
$$

% We single out the {\em Arf invariant}, given by
% $$
% \alpha = \sum_{i=1}^g \bar{a}_i \bar{b}_i
% $$
% as an example of a quadratic Boolean polynomial.
% Let $\Psi$ denote the subset
% of $ \Omega$ consisting of those Sp-forms with zero Arf invariant.
% In other words,
% $$
% \Psi = \{ \omega \in \Omega | \alpha(\omega) = 0 \}.
% $$

% Johnson \cite{Jo1} proved that Birman-Craggs homomorphisms are in
% bijective correspondence with elements of $\Psi$, and that
% $\Psi$ has $2^{2g-1} +
% 2^{g-1}$ elements.

Johnson proves in \cite{Jo1} that the image of $\sigma$ in the vector space of functions $\Omega
\rightarrow \ztwo$ is in fact isomorphic to $B_3$.  Thus we can write
$$\sigma: \Tgo \rightarrow B_3$$

\paragraph{Explicit formulas on generators.}
Recall that twists about separating curves generate $\kgo$ by definition.  In \S 7 of \cite{Jo1}, Johnson gives an explicit formula for $\sigma
(T_\gamma)$ when $\gamma$ is a separating curve on $\Sigma_{g,1}$.  The curve $\gamma$ separates $\Sigma_{g,1}$ into two components; let $\Sigma'$ denote the component which does not contain $\partial \Sigma_{g,1}$.
Let $A_1, B_1, \ldots, A_{g(\Sigma')},
B_{g(\Sigma')}$
be a symplectic $Z_2$-homology basis for $\Sigma'$.
Then Johnson's formula is:
\begin{equation}
\label{formula:sigma:2}
\sigma(T_{\gamma})=\sum_{i=1}^{g(\Sigma_\alpha)} \bar{A}_i
\bar{B}_i
\end{equation}
The expression on the right-hand side of the equation is
independent of the choice of symplectic basis for
$\Sigma_{\gamma\delta}$.  The restriction $\sigma|_{\kgo}$ is surjective onto $B_2 \subset B_3$.

Johnson gives a similar formula for generators of $\Tgo$.  We will not need this formula in our calculations, but we give it here for the sake of completeness.  Recall that a {\em bounding pair map} (BP map for short) in $\Tgo$ is the
composition $T_\gamma T_\delta^{-1}$ of Dehn twists about nonseparating
simple closed curves $\gamma,\delta$
whose union $\gamma\cup\delta$ separates $\Sigma_{g,1}$.  Given a BP map $T_\gamma T_\delta^{-1}$, we let $\Sigma'$ be the component of $\Sigma_{g,1} \setminus (\gamma \cup \delta)$ which does not contain $\partial \Sigma_{g,1}$, as above.
Let $C$ be the homology class of $\gamma$.  Then Johnson's formula for $\sigma:\Tgo\to B_3$ is as follows:
\begin{equation}
\label{formula:sigma:1}
\sigma(T_\gamma T_\delta^{-1})= \left(
\sum_{i=1}^{g(\Sigma')} \bar{A}_i \bar{B}_i \right)
(\bar{C} + 1)\in B_3
\end{equation}
Since BP maps generate $\Tgo$ (see \cite{Jo2}), this formula can be taken as the definition of the map $\sigma$.

\medskip
\noindent
{\bf Remark. }
As previously noted, $B_3$ is a $\mgo$-module in a fairly obvious way.  We
have  that $f \in Sp(2g, \Z_2)$ acts on a map $\omega:H \to \Z_2$
in $\Omega$ by $f \cdot \omega (x) = \omega(f(x))$.  Then $Sp$
acts on a function $\phi: \Omega \to \Z_2$ adjoint to its action
on $\Omega$, i.e., $f \cdot \phi (\omega) = \phi ( f \cdot
\omega)$.  Furthermore, $\sigma$ is a $\mgo$-equivariant map; that is,
for $f\in\Mod_{g,1}$ $h \in \Tgo$ we have
$\sigma (fhf^{-1}) = \hat{f} \cdot \sigma(h)$, where $\hat{f}$
denotes the image of the map $f$ under the symplectic
representation mod 2.
Moreover, since $\kgo$ is a $\mgo$-submodule of $\T$, we have that
$\sigma (\kgo) = B_2$ is also a $\mgo$-module and that
$\sigma|_{\kgo}$ is also a $\mgo$-equivariant map.

\medskip

Henceforth we will restrict our attention to the subgroup $\kgo$ in
$\Tgo$ and write simply $\sigma: \kgo \rightarrow B_2$ for the BCJ
homomorphism.  However, the tools and techniques we are about to discuss
will extend to the case of the full Torelli group.

\section{Abelian cycles in $H_2(\kgo,\ztwo)$}
\label{abeliancycles}

In this section, we will give a method for constructing nontrivial
classes in $H_2 (\kgo, \ztwo)$, inspired by Sakasai \cite{Sa}.
The idea is to construct abelian cycles in $H_2(\kgo,
\ztwo)$ and then to show that their images under the induced map
$$\sigma_*:H_2(\kgo,\ztwo)\to H_2(B_2, \ztwo) \cong \wedge^2 B_2 \oplus
B_2$$ are nontrivial.

\subsection{Abelian cycles}

If $f,g \in \kgo$ commute then we have a homomorphism
$i:\Z^2\to\kgo$, which induces a homomorphism
$$i_*:H_2(\Z^2,\ztwo)\to H_2 (\kgo, \ztwo)
 $$
The image of the generator $t$ of
$H_2(\Z^2, \ztwo)$ is denoted by
$$\{ f,g \}:= i_*(t) \in H_2(\kgo, \ztwo)$$
The homology class $\{ f,g \}$ is called the {\em abelian cycle}
corresponding to the pair $f,g$.

While we have given the definition of abelian cycles for
the case of $H_2(\kgo, \ztwo)$, one can generalize the construction of
abelian cycles in the obvious way for different choices of groups and of
coefficients, and also for degree $n > 2$ by taking $n$ commuting
elements and looking at the injection $\Z^n
\hookrightarrow \kgo$.

\paragraph{Homology of \boldmath$B_2$.}
Recall that as a vector space over $\F_2$, the summand $B_2$
of the algebra $B$ is spanned by $d=2g^2 + g + 1$ Boolean monomials.
Thus
$$H_2 (B_2, \ztwo) \cong \wedge^2 (B_2) \oplus B_2$$

The second summand comes from the Universal Coefficient Theorem and
the fact that $\textrm{Tor}(B_2, \ztwo) \cong B_2$.
Let $\{ f,g \}\in H_2(\kgo, \ztwo)$ be an abelian cycle.  Then it is
straightforward to verify (see e.g.\ Lemma~5.3 of \cite{Sa}) that
\begin{equation}
\label{split2}
\sigma_*(\{ f,g \}) = (\sigma(f) \wedge \sigma(g),0) \in \wedge^2
B_2 \oplus B_2
\end{equation}

Our strategy to prove that an abelian cycle in $H_2(\kgo,\Z)$ is nontrivial
will be to show that its image under $\sigma_\ast$ is nontrivial.

\bigskip
\noindent
{\bf Remark. }
It is impossible to detect classes corresponding to the $B_2$ summand of $H_2 (B_2, \ztwo)$ using this method, since an abelian
cycle arises from a map of $\Z^2$ into $\kgo$ and any such
class in $H_2 (\kgo, \ztwo)$ must vanish in the $\textrm{Tor}$ term
of the Universal Coefficient Theorem.
For this reason we will write $\wedge^2 B_2 \subset
H_2(B_2,\ztwo)$ and will refer to the element $(x \wedge y,0) \in
H_2(B_2,\ztwo) = \wedge^2 B_2 \oplus B_2$ simply as $x \wedge y
\in \wedge^2 B_2$.
Thus the formula (\ref{split}) becomes simply
\begin{equation}
\label{split}
\sigma_*(\{ f,g \}) = \sigma(f) \wedge \sigma(g) \in H_2(B_2, \F_2)
\end{equation}

\subsection{Calculations with \boldmath$\sigma_*$}
\label{mod2}

Our strategy for proving Theorem \ref{theorem:main1}
is twofold. First, we will use the induced map on
homology $\sigma_*$ together with the formula (\ref{split}) to hit basis
elements of $\wedge^2 B_2$ directly with abelian cycles. Secondly,
we will use the fact that $\sigma$ is a $\mgo$-equivariant map
(and hence so is $\sigma_*$) to increase our
efficiency. Thus we will compute orbits of basis elements of
$\wedge^2 B_2$ and note that hitting one element of an orbit with
$\sigma_*$ tells us that every basis element in the orbit is also
in the image of $\sigma_*$.

\paragraph{A basis for \boldmath$\wedge^2 B_2$.}

\begin{figure}[htpb!]
\begin{picture}(0,0)(0,11)
\put(130,70){ {\bf $a_1$} }
\put(171,70){ {\bf $a_2$} }
\put(211,70){ {\bf $a_3$} }
\put(275,70){ {\bf $a_g$} }
\put(130,30){ {\bf $b_1$} }
\put(171,30){ {\bf $b_2$} }
\put(213,30){ {\bf $b_3$} }
\put(276,30){ {\bf $b_g$} }
\end{picture}
\centerline{\includegraphics[scale=.5]{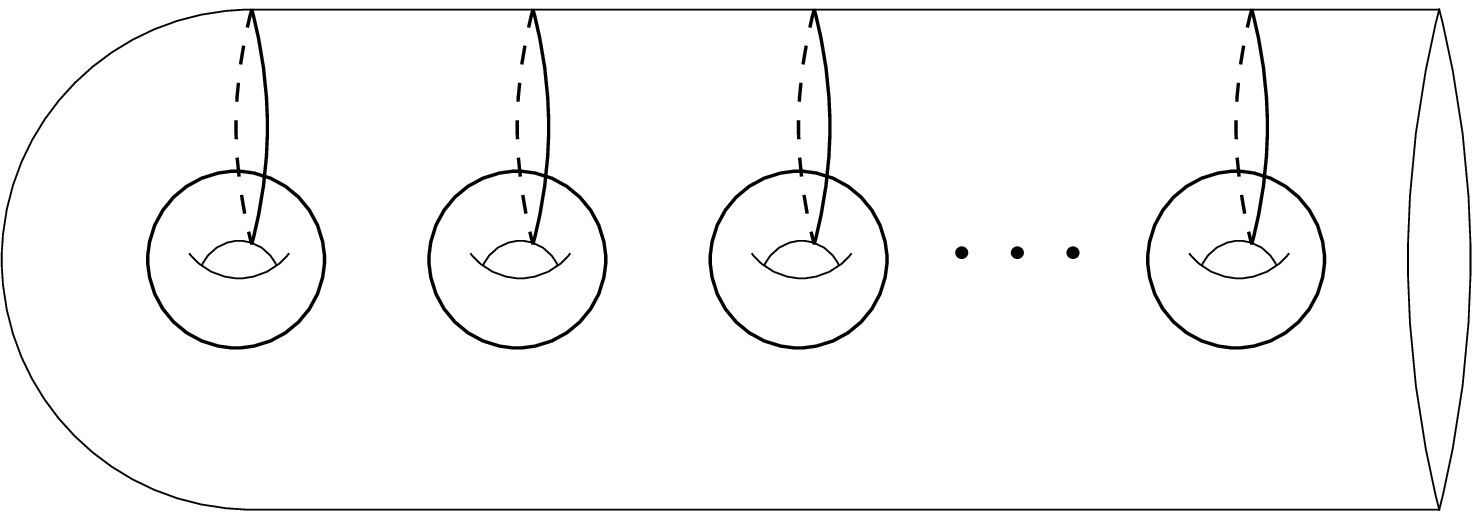}}
\caption{A fixed symplectic basis for $H_1 (S_{g,1}, \ztwo)$}
\label{symplecticbasis}
\end{figure}

Recall that a basis ${\cal B}$ for $B_2$ consists of the following $d =
2g^2 +g + 1$ Boolean monomials:
$${\cal B} = \{ 1,
\bar{a}_i, \bar{b}_i, \bar{a}_i \bar{b}_j, \bar{a}_i \bar{a}_j (i
\neq j), \bar{b}_i \bar{b}_j (i \neq j) \}$$ where the $a_i, b_j$
are the fixed symplectic basis for $H_1 (S_{g,1}, \ztwo)$ shown in Figure~\ref{symplecticbasis}.
Thus a
basis for $\wedge^2 B_2$ is given by
$$\left(
\begin{array}{c}
d \\
2
\end{array}
\right) = 2g^4 + 2g^3 +\frac{3}{2}g^2 +\frac{g}{2}$$ elements of the form $m_1
\wedge m_2$, where $m_1, m_2$ are distinct monomials in ${\cal
B}$.

When trying to describe a basis for $\wedge^2 B_2$, it is convenient to
organize elements according to distinct patterns of
the indices occurring in the monomials involved, which also have some topological significance.
We particularly wish to distinguish elements of the form
\begin{equation}
\bar{a}_i \bar{x} \wedge \bar{b}_i \bar{y}
\end{equation}
where $x$ and $y$ are arbitrary basis elements of $H$.  We refer to such a basis element of $\wedge^2 B_2$ as {\em
index-matched}
\label{index-matched}.  We allow that possibility that $x = a_i$ or $y =
b_i$, so that either or both factors of the wedge product could be
linear.

We are now able to state our main homology calculation.

\begin{theorem}
\label{image}
For $g \geq 4$, the image of the map $\sigma_*:H_2(\kgo,\F_2)\to
\wedge^2B_2$ contains the subspace $W$ spanned by all basis elements
which are not index-matched.  This space has dimension equal to a
polynomial $p(g) =\lowerbound$.
\end{theorem}

\proof
We want to examine the image of the $\Mod_{g,1}$-equivariant map
$\sigma_*$ in $B_2$.  One approach would be to calculate equivalence
classes of basis elements of $B_2$ under the full action of
$\Mod_{g,1}$.
However, it turns out to be more efficient to consider only the action
of two particularly simple $\Mod_{g,1}$-maps on the $\F_2$-vector space
$B_2$:
\begin{equation}
\label{eq:map1}
\bar{a}_i \mapsto \bar{b}_i \mapsto \bar{a}_i
\end{equation}
\begin{equation}
\label{eq:map2}
\bar{a}_i \mapsto \bar{a}_j; \hskip .15in \bar{b}_i \mapsto \bar{b}_j
\end{equation}

Elements not mentioned in the above maps are understood to be fixed.
Recall that $W$ is the subspace of $B_2$ spanned by all basis elements
which are not index-matched.  What we do now is to
calculate equivalence classes of
basis elements of $W$ under the two maps given by (\ref{eq:map1}) and
(\ref{eq:map2}).  There are $11$
equivalence classes of basis elements of $W$.  We
refer to these $11$ equivalence
classes as {\em partial Sp-orbits}.  They are:

\vskip .15in

\noindent
For all $i,j,k,l = 1, \ldots, g, \hskip .1in  i \neq j \neq k \neq l$:
\begin{eqnarray*}
\textrm{I} &=& \{\bar{a}_i \bar{b}_i \wedge \bar{a}_j \bar{b}_j \} \\
\textrm{II} &=& \{\bar{a}_i \bar{b}_i \wedge \bar{a}_j \bar{b}_k, \bar{a}_i \bar{b}_i \wedge \bar{a}_j \bar{a}_k, \bar{a}_i \bar{b}_i \wedge \bar{b}_j \bar{b}_k \} \\
\textrm{III} &=& \{ \bar{a}_i \bar{a}_j \wedge \bar{a}_k \bar{a}_l, \bar{a}_i \bar{a}_j \wedge \bar{a}_k \bar{b}_l, \bar{a}_i \bar{a}_j \wedge \bar{b}_k \bar{b}_l, \bar{a}_i \bar{b}_j \wedge \bar{a}_k \bar{b}_l, \bar{a}_i \bar{b}_j \wedge \bar{a}_k \bar{b}_l, \bar{b}_i \bar{b}_j \wedge \bar{b}_k \bar{b}_l \}\\
\textrm{IV} &=&  \{\bar{a}_i \bar{a}_j \wedge \bar{a}_i \bar{a}_k, \bar{a}_i \bar{a}_j \wedge \bar{a}_i \bar{b}_k, \bar{a}_i \bar{b}_j \wedge \bar{a}_i \bar{b}_k, \bar{a}_i \bar{b}_j \wedge \bar{a}_k \bar{b}_j, \bar{a}_i \bar{b}_j \wedge \bar{b}_j \bar{b}_k, \bar{b}_i \bar{b}_j \wedge \bar{b}_i \bar{b}_k\}\\
\textrm{V} &=&  \{ \bar{a}_i \wedge \bar{a}_j \bar{b}_j, \bar{b}_i \wedge \bar{a}_j \bar{b}_j \} \\
\textrm{VI}&=& \{ \bar{a}_i \wedge \bar{a}_i \bar{a}_j, \bar{a}_i \wedge \bar{a}_i \bar{b}_j, \bar{b}_i \wedge \bar{a}_j \bar{b}_i, \bar{b}_i \wedge \bar{b}_i \bar{b}_j \}   \\
\textrm{VII}&=& \{ \bar{a}_i \wedge \bar{a}_j \bar{a}_k, \bar{a}_i \wedge \bar{a}_j \bar{b}_k, \bar{a}_i \wedge \bar{b}_j \bar{b}_k, \bar{b}_i \wedge \bar{a}_j \bar{a}_k, \bar{b}_i \wedge \bar{a}_j \bar{b}_k, \bar{b}_i \wedge \bar{b}_j \bar{b}_k \} \\
\textrm{VIII}&=&  \{ 1 \wedge \bar{a}_i \bar{b}_i \}     \\
\textrm{IX}&=&  \{ 1 \wedge \bar{a}_i \bar{a}_j, 1 \wedge \bar{a}_i \bar{b}_j, 1 \wedge \bar{b}_i \bar{b}_j \}    \\
%\textrm{OLDLL1}&=& \{ \bar{a}_i \wedge \bar{a}_i, \bar{b}_i \wedge \bar{b}_i \}     \\
\textrm{X}&=& \{ \bar{a}_i \wedge \bar{a}_j, \bar{a}_i \wedge \bar{b}_j, \bar{b}_i \wedge \bar{b}_j \}\\
\textrm{XI}&=& \{ 1 \wedge \bar{a}_i, 1 \wedge \bar{b}_i \}
\end{eqnarray*}

\paragraph{Direct calculations.}
We will work our way through the above-listed orbits one at a time,
building on our previous calculations as we go.  We will give a few
calculations explicitly in order to give the reader the idea of how to
proceed.  For the remaining cases, we will give the readers the
required abelian cycles and leave the calculations as an exercise.

We begin by noting a fact which will simplify our calculations.
Let $\alpha, \beta$ be two simple closed curves on a surface such that
$i(\alpha, \beta) = 1$.  Then the boundary $\gamma$ of a regular
neighborhood $N$ of $\alpha \cup \beta$ is a genus $1$ separating curve.
Further, the two curves $\alpha$ and $\beta$ form a symplectic basis for
$H_1(N, \ztwo)$, and hence $\sigma(T_\gamma) =\bar{\alpha} \bar{\beta}$
(in general we will not distinguish between a curve and its homology
class).

We call the pair $\alpha, \beta$ a {\em spine} for the genus 1
separating curve $\gamma$.  We will also abuse terminology and refer to
the spine of the map $T_\gamma$ and the spine of the subsurface bounded
by $\gamma$.  Thus when trying to ``hit'' a specific basis element of
$\wedge^2 B_2$, we will in general look first for spines which give the
desired monomials.  Clearly, disjoint spines
correspond to disjoint separating curves and hence to commuting twists
in $\kgo$.

\begin{figure}[htpb!]
\begin{picture}(0,0)(0,11)
\put(180,71){ {\bf $a_i$} }
\put(217,71){ {\bf $a_j$} }
\put(173,30){ {\bf $b_i$} }
\put(210,30){ {\bf $b_j$} }
\end{picture}
\centerline{\includegraphics[scale=.5]{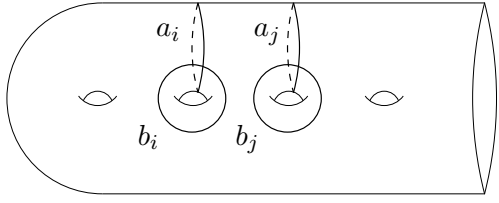}}
\caption{The spines of two bounding curves corresponding to an abelian
cycle in $H_2(\kgo, \ztwo)$ which maps to an element of Partial orbit I.}
\label{figure:qq1}
\end{figure}

\paragraph{Partial orbit I.}
Referring to Figure~\ref{figure:qq1}, let $\gamma_i, \gamma_j$ be the
two separating curves corresponding to the two `spines' shown on the
$i^{th}$ and $j^{th}$ holes, respectively.   Then $\{ T_{\gamma_i},
T_{\gamma_j} \}$ is an abelian cycle in $H_2(\kgo, \ztwo)$, and using
(\ref{formula:sigma:2}) and (\ref{split}) gives
\begin{eqnarray*}
\sigma_2 ( \{ T_{\gamma_i}, T_{\gamma_j} \} ) &=&
\sigma(T_{\gamma_i}) \wedge \sigma(T_{\gamma_j}) \\
&=& \bar{a}_i \bar{b}_i \wedge \bar{a}_j \bar{b}_j
\end{eqnarray*}

\begin{figure}[htpb!]
\begin{picture}(0,0)(0,11)
\put(182,70){ {\bf $a_i$} }
\put(219,70){ {\bf $a_j$} }
\put(175,28){ {\bf $b_i$} }
\put(212,30){ {\bf $b_j + a_k$} }
\end{picture}
\centerline{\includegraphics[scale=.5]{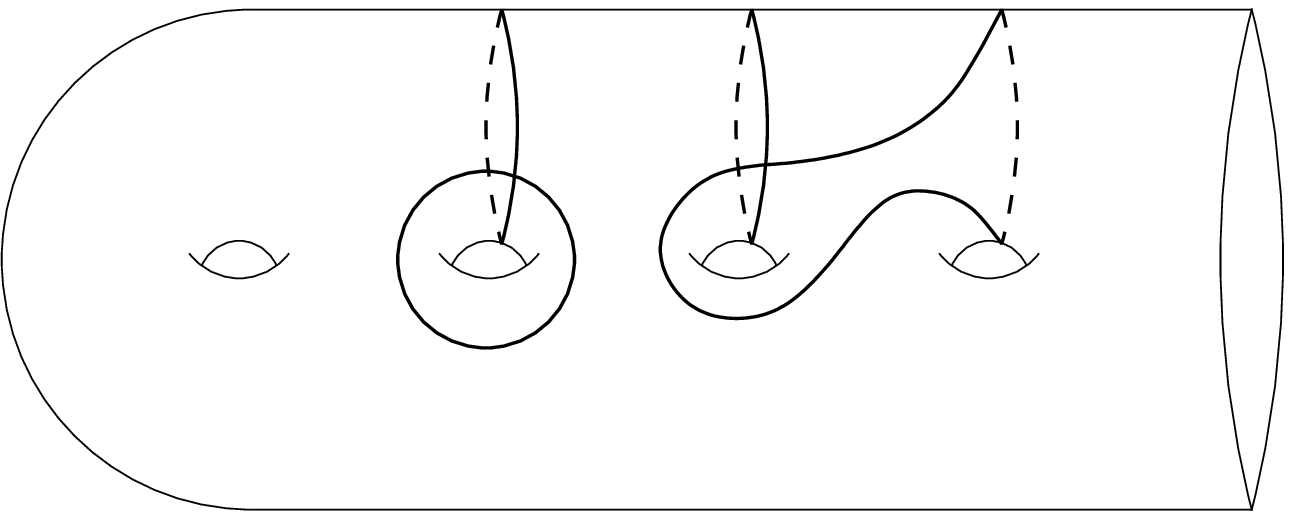}}
\caption{The spines of two bounding curves corresponding to an abelian
cycle mapping to Partial orbit II.}
\label{figure:qq2}
\end{figure}

\paragraph{Partial orbit II.}  Let $\gamma, \delta$
denote the two separating curves corresponding to the spines shown in
Figure~\ref{figure:qq2}; let the holes being `used' in the picture be
the $i^{th}, j^{th}$, and $k^{th}$ holes, respectively.  Then $\{
T_{\gamma}, T_\delta \}$ is an abelian cycle, and
\begin{eqnarray*}
\sigma_2 ( \{ T_{\gamma}, T_{\delta} \} ) &=&
\sigma(T_{\gamma}) \wedge \sigma(T_{\delta}) \\
&=& \bar{a}_i \bar{b}_i \wedge \bar{a}_j \bar{c} \\
&=& \bar{a}_i \bar{b}_i \wedge \bar{a}_j (\bar{b_j + a_k}) \\
&=& \bar{a}_i \bar{b}_i \wedge \bar{a}_j (\bar{b}_j + \bar{a}_k + b_j \cdot a_k) \\
&=& \bar{a}_i \bar{b}_i \wedge (\bar{a}_j \bar{b}_j + \bar{a}_j \bar{a}_k) \\
&=& \bar{a}_i \bar{b}_i \wedge \bar{a}_j \bar{b}_j + \bar{a}_i
\bar{b}_i  \wedge \bar{a}_j \bar{a}_k
\end{eqnarray*}
Since the first term on the right-hand side of the last equality is
already in the image of $\sigma_2$ (Partial orbit I), so is the second term.

\begin{remark}
Note that there exists $h\in \Mod_{g,1}$ taking the spine (hence the
corresponding abelian cycle) illustrated in
Figure~\ref{figure:qq1} to that illustrated in Figure~\ref{figure:qq2}.
 Hence these abelian cycles lie in the same $\Sp$-orbit, namely
(partial) orbit I.
The calculation just made exhibits an element from (partial) orbit II as a
difference of these two cycles.
\end{remark}

\smallskip
We next give an example showing how to hit a linear term.

\paragraph{Partial orbit V.}  \begin{figure}[htpb!]
\begin{picture}(0,0)(0,11)
\put(169,72){ {\bf $a_i + b_i$} }
\put(182,55){ {\bf $a_i + a_j$} }
\put(240,70){ {\bf $a_k$} }
\put(247,45){ {\bf $b_k$} }
\end{picture}
\centerline{\includegraphics[scale=.5]{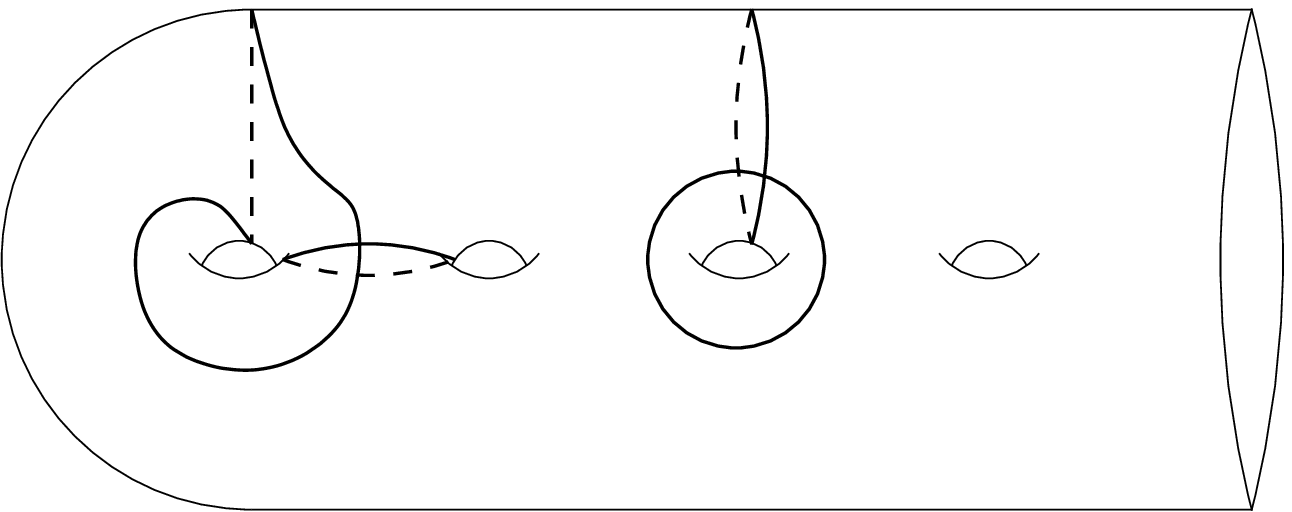}}
\caption{The spines of two bounding curves corresponding to an abelian
cycle mapping to partial orbit V.}
\label{figure:lq1}
\end{figure}

Let $\gamma, \delta$ denote the two separating curves corresponding to
the spines shown in Figure~\ref{figure:lq1}; let the holes being `used'
in the picture be the $i^{th}, j^{th}$, and $k^{th}$ holes,
respectively.  Then $\{ T_{\gamma}, T_\delta \}$ is an abelian cycle,
and
\begin{eqnarray*}
\sigma_2 ( \{ T_{\gamma}, T_{\delta} \} ) &=&
\sigma(T_{\gamma}) \wedge \sigma(T_{\delta}) \\
&=& (\bar{a_i + b_i}) (\bar{a_i + a_j}) \wedge \bar{a}_k \bar{b}_k \\
&=& (\bar{a}_i + \bar{b}_i + 1)(\bar{a}_i + \bar{a}_j) \wedge \bar{a}_k
\bar{b}_k \\ &=& (\bar{a}_i \bar{b}_i + \bar{a}_i \bar{a}_j + \bar{a}_j
\bar{b}_i + \bar{a}_j)\wedge \bar{a}_k \bar{b}_k \\ &=& \bar{a}_i
\bar{b}_i \wedge \bar{a}_k \bar{b}_k + \bar{a}_i \bar{a}_j \wedge
\bar{a}_k \bar{b}_k + \bar{a}_j \bar{b}_i \wedge \bar{a}_k \bar{b}_k +
\bar{a}_j \wedge \bar{a}_k \bar{b}_k
\end{eqnarray*}

Since the first three terms on the right-hand side of the last equality
come from partial orbits I and II, it follows that partial
orbit V is also contained in
the image of $\sigma_*$.

The remaining calculations required for the
proof of Theorem~\ref{image} are similar in flavor, although
lengthier; we defer them to the Appendix, where we
present figures
indicating the necessary abelian cycles corresponding to
each of the remaining (partial) orbits listed above.

\begin{remark}  In many of our calculations, such as the three given above,
we do not use the hypothesis that $g \geq 4$; lower bounds may thus be
obtained for $g<4$.
\end{remark}

\bigskip
\noindent
{\bf Cokernel.}  Let $\IM$ denote the subspace of $\wedge^2 B_2$ spanned
by index-matched basis elements.  If one considers an index-matched
element of the form $\bar{a}_i x \wedge \bar{b}_i y$ where $x$ and $y$
are distinct basis element with distinct indices, then it is clear that
the dimension of $\IM$ is a polynomial in $g$ with highest degree term
$4g^3$.  However, it turns out that ``most'' of $\IM$ is actually
contained in the image of $\sigma_*$.

\begin{proposition}
\label{prop:cokernel}
The cokernel of $\sigma_*:H_2(\kgo,\F_2)\to H_2(B_2,\F_2)$
has dimension at most a polynomial with highest degree term $4g^2$.
\end{proposition}

\proof
Theorem~\ref{image} above tells
us that the dimension of the cokernel of $\sigma_*$
is at most $\dim(\IM)$.  Index-matched elements of the
type $\bar{a}_i x \wedge \bar{b}_i y$, where $x$ and $y$ are distinct
basis element with distinct indices, are the only index-matched basis
elements contributing a cubic term to $\dim(\IM)$.
The subspace spanned by such elements has dimension
$g(2g-2)(2g-3)$

Figure~\ref{figure:im1} below shows two abelian cycles. The first, when
taken together with Theorem~\ref{image}, shows that the image of
$\sigma_*$ contains any sum of the form
\begin{equation}
\label{eq:lastcalc}
\bar{a}_i \bar{b}_i
\wedge \bar{a}_i \bar{b}_j + \bar{a}_j \bar{b}_j \wedge \bar{a}_i \bar{b}_j
\end{equation}

The second, when combined with Theorem~\ref{image} and
(\ref{eq:lastcalc}), and after an adjustment of indices, shows that
any sum of the form
$$
a_i b_j \wedge b_i a_k + a_l b_j \wedge b_l a_k
$$
(where all distinct indices are assumed to be not equal) is contained in the image of $\sigma_*$.
The subspace spanned by this element and
its $\Mod_{g,1}$-orbit has dimension $(g - 1)(2g - 2)(2g -3)$.  It
follows that index-matched elements of the above type
contribute at most $4g^2-10g+6$ to the dimension of the
cokernel of $\sigma_*$.
\begin{figure}[htpb!]
\begin{picture}(0,0)(0,11)
%\put(162,70){ {\bf $a_i$} }
%\put(199,70){ {\bf $a_j$} }
%\put(155,25){ {\bf $b_i$} }
%\put(192,30){ {\bf $b_j + a_k$} }
\end{picture}
\centerline{\includegraphics[scale=1]{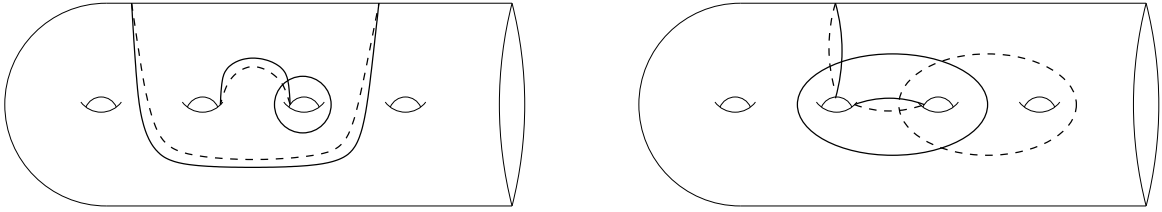}}
\caption{Two spines giving sums of pairs of index-matched elements.}
\label{figure:im1}
\end{figure}

Next, we consider index-matched basis elements
whose orbits contribute a $g^2$-term to the dimension of $
\IM$.  Calculations similar to the above show that the image
of $\sigma_*$ also contains sums of these basis elements spanning a
subspace whose dimension has the same quadratic term.  Thus the
remaining index-matched basis elements contribute at most a linear term
to the dimension of the cokernel of $\sigma_*$.  This completes the
proof of the
proposition.\endproof

\section{Integral abelian cycles and the Casson-Morita algebra}

We next investigate the question of whether the mod 2 abelian cycles
constructed in the previous section can be lifted to integral classes.

\subsection{The Casson-Morita algebra}

As described in \S\ref{spforms} above, the $\F_2$-algebra $B_3$
is a kind of ``universal
receptor'' for Rochlin invariants, or more precisely
Birman-Craggs homomorphisms (see \S\ref{spforms} above).  We now
explain how Morita lifted this setup to $\Z$ using the Casson invariant.

\paragraph{Casson invariant.}
The Casson invariant $\lambda(M)$ of a homology 3-sphere $M$ is, roughly
speaking, half the algebraic number of conjugacy classes of irreducible
representations of $\pi_1(M)$ into SU(2). The invariant $\lambda$ is
integer-valued, and is a lifting of Rochlin's $\ztwo$-valued invariant,
defined in Section~\ref{spforms}, in the sense that $\lambda(M) \equiv
\mu(M) \hskip .05in \textrm{mod } 2$.  See \cite{AM} for a
thorough exposition of the Casson invariant.

Let
$h:\Sigma_g \rightarrow S^3$
be a Heegaard embedding, and let $f \in \T$.  Now split $S^3$
along $h(\Sigma_g)$ and reglue via the map $f$.  By Meyer-Vietoris,
the resulting
3-manifold $M(h,f)$ is a homology 3-sphere.  Just as
described in Section \ref{spforms} above for the Rochlin
invariant, one can define a function  $\lambda_h: \T\rightarrow \Z$
given by $f\mapsto \lambda(M(h,f))$.
%\begin{eqnarray*}
%\lambda_h: \T &\rightarrow& \Z\\
%f &\mapsto& \lambda(M(h,f))
%\end{eqnarray*}

While $\lambda_h$ is generally not a homomorphism, Morita \cite{Mo1,Mo2}
proved that the restriction $\lambda_h: \kg \rightarrow \Z
$ is a homomorphism.  He also showed that every homology 3-sphere may be
represented as $M_\phi$ for some $\phi \in \kg$.  See \cite{Mo1}.

\paragraph{The Casson-Morita algebra. }

The method of combining all the homomorphisms $\lambda_h$ into
one package is motivated by Casson's knot invariant.
Let $K$ be a knot in a homology sphere
$M$, and let $M_K$ denote the homology 3-sphere which results from doing
$(1,1)$-surgery on $K$.  Define $\lambda'(K):=\lambda(M_K) -
\lambda(M)$.
% Casson proved that the knot invariant $\lambda'$ satisfies
% $\lambda'(K)=\frac{1}{2} \Delta''_K(1)$,
% where $\Delta''_K$ denotes the second derivative
% of the (symmetrized) Alexander polynomial of the knot $K$.
Morita \cite{Mo1,Mo2}
has shown that if $c$ is a separating curve on $\Sigma_{g,1}$ then
$$\lambda_h(T_c) = - \lambda'(h(c))$$

Now consider the Seifert matrix $L = (l_{ij})$ of the knot $K$, given
with respect to some choice of symplectic basis for the homology of a
Seifert surface for $K$.  By definition, $l_{ij}$ is the linking number
of the $i^{th}$ basis element with the positive push-off of the $j^{th}$
basis element. It turns out that $\lambda'(K)$ can be expressed as a
polynomial of degree 2 in the $l_{ij}$; this polynomial is
given explicitly in \cite{Mo1}.

The {\em Casson-Morita algebra} $\calA$ is then defined as
follows.  Let $\calA$ be the commutative $\Z$-algebra generated by the
abstract symbols $l(u,v)$, where $u,v \in H_1(S)$, satisfying the
relations:
\begin{enumerate}
\item $l(v,u) = l(u,v) + u \cdot v$
\item $l(n_1 u_1 + n_2 u_2, v) = n_1 l(u_1, v) + n_2 l(u_2, v)$ \
\mbox{for $m,n\in \Z$}
\end{enumerate}

For each Heegaard embedding $h: S \rightarrow S^3$, the algebra $\calA$ has an
{\em evaluation homomorphism} $\epsilon_h: \calA \rightarrow \Z$ given
by $$\epsilon_h(l(u,v))=lk(h_*(u),h_*(v)^+)$$

\subsection{Morita's homomorphism \boldmath$\rho$}

The {\em Morita homomorphism} $\rho: \kgo \rightarrow \calA$,
introduced in \cite{Mo1}, is defined as follows.  Given a generator $T_c
\in \kgo$, such that the curve $c$ bounds a subsurface $S'$ with
symplectic homology basis $A_1, \ldots, A_g, B_1, \ldots, B_g$, we set
\begin{eqnarray}
\rho(T_c) = &-& \sum_{i=1}^{g(S')} [ l(A_i, A_i) l(B_i, B_i) -
l(A_i,B_i)l(B_i,A_i)] \nonumber \\ &-& 2 \sum_{i<j\leq g(S')} [ l(A_i,
A_j) l(B_i, B_j) - l(A_i,B_j)l(A_j,B_i)]
\label{rhoformula}
\end{eqnarray}

Morita \cite{Mo1} proved this assignment extends to a homomorphism on
all of $\kgo$.  The map $\rho$ plays a similar role
with regard to the Casson invariant to that of the BCJ map with regard
to the Rochlin invariant; we will describe the relationship of $\rho$
with the Casson
invariant explicitly in Proposition~\ref{bigdiagram} below.  In fact,
Morita's homomorphism $\rho$ is actually a lift of the BCJ homomorphism,
just as the the Casson invariant is a lift of the Rochlin invariant.

\begin{proposition}
\label{triangle}
There exists a homomorphism
$\mu: \A \rightarrow B_2$ such that the following diagram commutes:
$$\xymatrix{
& \A \ar[d]^{\mu} \\
\kgo \ar[ur]^{\rho} \ar[r]^{\sigma} & B_2  }
$$
\end{proposition}

\medskip
We would like to point out that the ``reduction homomorphism'' $\mu$ must respect the
relations in the algebra ${\cal
A}$, and so takes a more complicated form than a simple ``reduction mod
$2$''.

\proof
Before making any assignments for the value of $\mu$ on any of the
generating symbols $l(a,b)$, we make a few motivational observations.
First, note that the image of the second term of
Equation~\ref{rhoformula} under any homomorphism to a $2$-group
must necessarily be 0.
Also, we claim that the first defining relation in $\A$ gives that the
image of the second part of the first term of Equation~\ref{rhoformula}
must also map to 0, since we
have:
\begin{eqnarray*}
l(A_i, B_i) l(B_i, A_i) &=& l(A_i, B_i) [l(A_i,B_i) + A_i \cdot B_i] \\
&=&l(A_i, B_i)^2 +  (A_i \cdot B_i) l(A_i, B_i)\\
&=& l(A_i, B_i)^2 + l(A_i, B_i)
\end{eqnarray*}
The square-free condition in $B_2$ now implies the claim.

Thus we can focus our attention on the image of $l(A_i,A_i)$.  Now, the obvious assignment which makes the diagram commutative is
$$
\mu(l(A_i, A_i)) = \bar{A}_i \in B_2
$$
The question becomes: how does one define the map $\mu$ on the entire algebra $\A$ so as to achieve the desired result?

Let $a_1, \ldots, a_g, b_1, \ldots, b_g$ be the fixed symplectic basis
for $H_1(S)$.  Then the defining relations show that symbols of the form
$l(x,y)$, where $x,y \in\{a_1, \ldots, a_g, b_1, \ldots, b_g \}$,
suffice to generate $\A$.  We now set, for all $1\leq i,j\leq g$:
$$\begin{array}{l}
\mu(l(a_i, a_i))=\bar{a}_i, \ \ \ \mu(l(b_i, b_i))=\bar{b}_i \\
\mu(l(a_i, b_j))=0 \\
\mu(l(b_j, a_i))=\mu(l(a_i, a_j))=\mu(l(b_i, b_j))=\delta_{ij}
\end{array}$$
A straightforward check shows
that these assignments satisfy the defining relations of $\A$, and thus
by extending linearly we get a well-defined homomorphism $\mu:{\cal
A}\to B_2$.  Further, it is also straightforward
(although a bit tedious) to check
that this definition on the fixed basis also gives $ \mu(l(u,u)) = \bar{u}
$ for any $u \in H_1(S)$, and thus the diagram commutes.
\endproof

We are now ready to prove the following.

\begin{proposition}\label{bigdiagram}The following diagram commutes:\newline
\xymatrix{
\hskip 1in &  \kgo \ar[0,2]^{\rho} \ar[d]_{\pi} \ar@/^2pc/[0,4]^\sigma &
  & \A \ar[d]^{\epsilon_f} \ar[0,2]^{\mu} &
  & B_2 \ar[d]^{\bar{\epsilon}_f} \\
\hskip 1in &  \kg \ar[0,2]^{\lambda_f} &
  & \Z \ar[0,2]^{\textrm{mod 2}} &
  & \ztwo      }
\end{proposition}

\proof
The fact that the left-hand square commutes is Theorem 2.2 of
\cite{Mo2}.  The
commutativity of the upper `triangle' is Proposition~\ref{triangle}
above. It remains to deal with the right-hand square.

Let $\omega_f$ denote the mod 2 self-linking form on $H_1(S)$ induced by
the Heegaard embedding $f: \Sigma_{g,1} \rightarrow M$.  Recalling the definition of
elements of $B_2$ from Section~\ref{spforms}, we have the following
natural analog of an evaluation map on $B_2$:
\begin{eqnarray*}
\bar{\epsilon}_f(\bar{u}) &=& \bar{u}(\omega_f) \\
&=& lk(f_*(u), f_*(u)^+) \mbox{ mod 2}
\end{eqnarray*}
It is now straightforward to check the commutativity of the right-hand square; again, it suffices to check on the generators $l(a_i, b_j)$, etc.  This finishes the proof of the proposition. \endproof

\subsection{Lifting abelian cycles}

We next describe how the $\ztwo$-classes constructed as in
Section~\ref{mod2} can be lifted to integral classes.  Define
$p:H_2( \kgo, \Z) \rightarrow H_2(\kgo, \ztwo)$ to be the composition
$$H_2( \kgo, \Z) \rightarrow H_2( \kgo, \Z) \otimes \ztwo
\rightarrow H_2(\kgo, \ztwo)$$
where the first map is given by $f \mapsto f \otimes 1$ and the second
map is the injection given by the
Universal Coefficients Theorem.  The analogous map $q: H_2(B_2, \Z)
\rightarrow H_2(B_2, \ztwo)$ is just the injection arising from the
Universal Coefficients Theorem.  Proposition~\ref{bigdiagram} directly
implies the following.

\begin{proposition}
The following diagram commutes.

\

\xymatrix{
\hskip .5in &  H_2( \kgo, \Z) \ar[r]^{\rho_*} \ar[d]_p & H_2(\A, \Z) \ar[r]^{\mu_*} & H_2 (B_2, \Z) \ar[d]^q \\
\hskip .5in &  H_2(\kgo, \ztwo) \ar[0,2]^{\sigma_*} & & H_2(B_2, \ztwo) }
\end{proposition}

Note that, by construction, the map $p$ sends an abelian cycle $\{f,
g\}$ in $H_2( \kgo, \Z)$ to the `same' abelian cycle in $H_2(\kgo,
\ztwo)$.  The proposition tells us that if $\xi \in H_2( \kgo, \Z)$ and
$\sigma_* p(\xi) \neq 0$, then $\xi \neq 0$ as well.  In other words,
each of the $\ztwo$ classes arising from abelian cycles which were
constructed in Section~\ref{mod2} and in the Appendix also lift to
integral classes.
Thus we also obtain lower bounds on the rank of
$H_2(\kgo, \Z)$ as a corollary to Theorem~\ref{image}.

\begin{corollary}
The rank of $H_2(\kgo, \Z)$ is at least $\lowerbound$.
\end{corollary}

\begin{remark}
In light of the relationship between the maps $\rho$ and $\sigma$, it is
natural to ask whether we could gain any information about $H_2(\kgo,
\Z)$ by working directly with integer coefficients, rather than with
$\F_2$.  However, note that
the $\Z$-rank of $H_2(\A, \Z)$ is precisely equal to the
$\ztwo$-dimension of $H_2(B_2, \ztwo)$.  Thus we lost no information by
working with $\sigma$ and with $\ztwo$ coefficients.  Morever, the
calculations required for the proof of Theorem~\ref{image} were greatly
simplified by working mod 2.
\end{remark}

\section{Cohomology}

The following commutative diagram summarizes the various relationships
between $H_*(\kgo, \ztwo)$ and $H^*( \kgo, \ztwo)$, the induced maps
$\sigma_*$ and $\sigma^*$, and their images.

\xymatrix{
\hskip .5in & H^2( B_2, \ztwo) \ar[0,2]^{\sigma^*} \ar[2,0]_{\cong}& &H^2 (\kgo, \ztwo) \ar[2,0]^{\cong}\\
& & & \\
\hskip .5in &  (H_2(B_2, \ztwo))^\star \ar[0,2]^{(\sigma_*)^\star} & &(H_2(\kgo, \ztwo))^\star \\
& & & \\
\hskip .5in & H_2(B_2, \ztwo) \ar[-2,0]^{\cong}& &H_2(\kgo,\ztwo) \ar[0,-2]_{\sigma_*}\ar@{-->}[-2,0]
}

\bigskip
Here $\sigma_*^\star$ denotes the dual of the
induced map $\sigma_*$ on homology.  The commutativity
of the top square comes directly from the Universal Coefficient Theorem
for cohomology.  We emphasize that we do not
necessarily get an isomorphism between $H_2(\kgo,\ztwo)$ and $(H_2(\kgo,
\ztwo))^\star$ since it is not known whether or not $H_2(\kgo, \ztwo)$
is finite dimensional.

Given $\xi_* \in H_2(B_2, \ztwo)$, we let $\xi^*$ denote its isomorphic
image in $H^2(B_2, \ztwo)$.  It now follows easily from the diagram that $$
\xi_* \in \mbox{ Im $\sigma_*$} \hskip .15in \Rightarrow \hskip .15in \xi^* \notin \mbox{ ker $\sigma^*$}
$$
Thus Theorem \ref{image} implies Theorem
\ref{theorem:main1}.     In fact, one can determine from the
nontriviality of the abelian cycle $\{f,g\}\in H_2(B_2,\F_2)$ the
nontriviality of the image under $\sigma^\ast$ of a
corresponding cup product of elements in $H^1(B_2,\F_2)$.

\bigskip
\noindent
{\bf Further classes. }
As previously noted, all abelian cycles constructed in $H_2(\kgo,\F_2)$
survive in $H_2( \T, \ztwo)$.  Of course,
one could also construct more classes by
utilizing Johnson's formula for bounding pair maps in addition to his
formula for twists about separating curves.  This combinatorial
challenge has been successfully undertaken by Vijay Ravikumar \cite{Ra}, who is
able to prove the following:

\begin{proposition}[Ravikumar]
The image under $\sigma_*$ of the full Torelli group contains a subspace of dimension $64g^6 + \mbox{lower-order terms}$.
\end{proposition}

This $g^6$ lower bound is the best order-of-magnitude
one can hope for using the method of abelian cycles.

\renewcommand{\theequation}{A-\arabic{equation}}
  % redefine the command that creates the equation no.
\setcounter{equation}{0}  % reset counter
\section*{Appendix: Calculations for the proof of Theorem~\ref{image}}\label{appendix}  % use *-form to suppress numbering

In this appendix we list figures indicating the abelian cycles
necessary to complete the proof of Theorem~\ref{image}. In each case,
the pictures show two spines corresponding to an abelian
cycle in $H_2(\kg, \ztwo)$, as in the proof of Theorem~\ref{image}
above.

We also remind the reader that a given figure labelled by a
particular Sp-orbit may not give an element of that orbit exactly,
but rather may correspond to a sum of an element in that orbit
together with basis elements already known to be in the image of
$\sigma_*$.

\vskip .2in

\begin{figure}[htpb!]
\begin{picture}(0,0)(0,11)
\put(40,300){ {\large \bf Orbit I} }
\put(195,300){ {\large \bf Orbit II} }
\put(345,300){ {\large \bf Orbit III} }
\put(40,200){ {\large \bf Orbit IV} }
\put(195,200){ {\large \bf Orbit V} }
\put(345,200){ {\large \bf Orbit VI} }
\put(40,100){ {\large \bf Orbit VII} }
\put(195,100){ {\large \bf Orbit VII} }
\put(345,100){ {\large \bf Orbit IX} }
\put(110,-5){ {\large \bf Orbit X} }
\put(290,-5){ {\large \bf Orbit XI} }
\end{picture}
\centerline{\includegraphics[scale=.47]{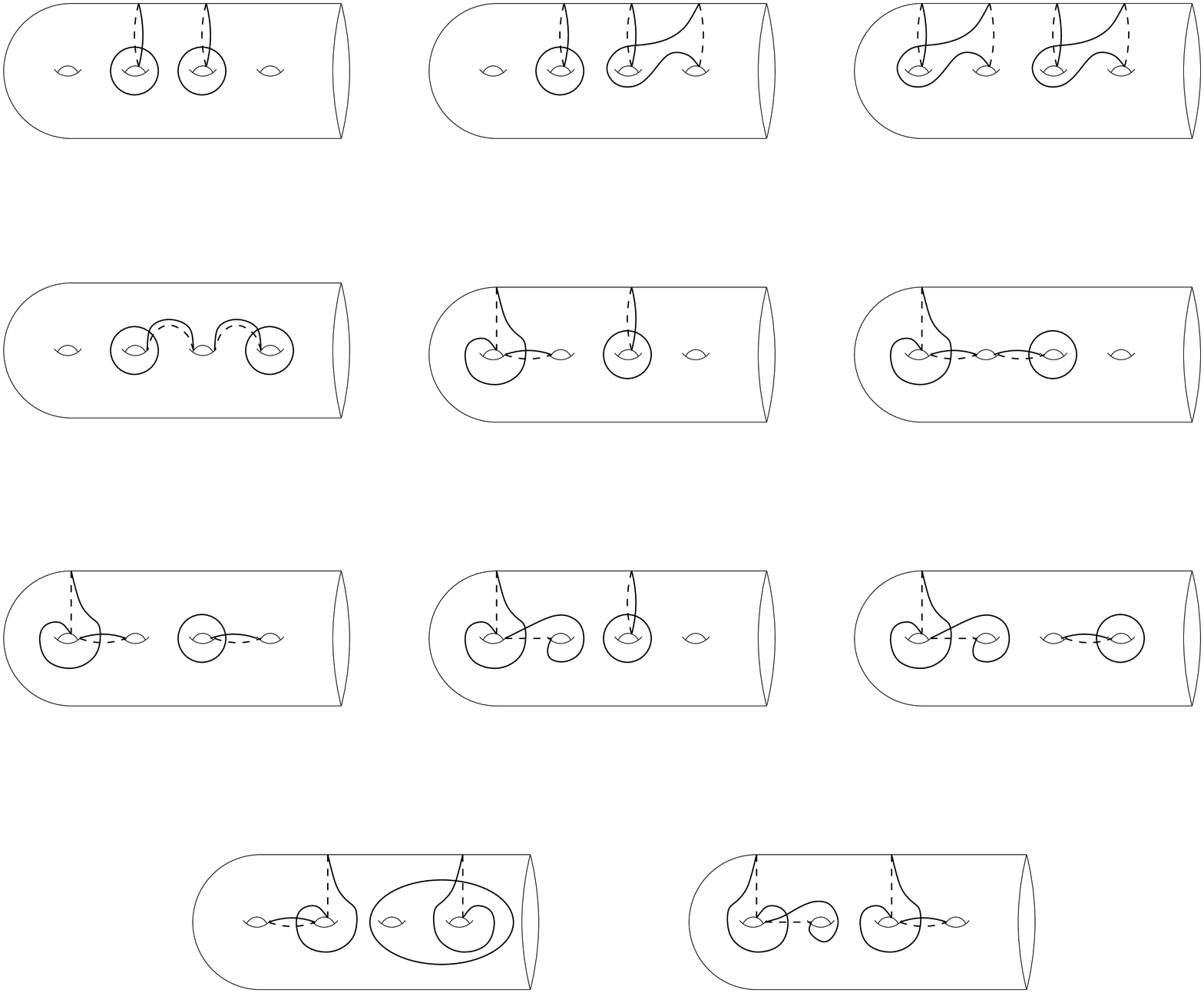}}
\label{figure:list1}
\end{figure}

% \newpage
% \vskip 1in

% \begin{figure}[htpb!]
% \begin{picture}(0,0)(0,11)
% \put(65,235){ {\large \bf Orbit VII} }
% \put(265,235){ {\large \bf Orbit VIII} }
% \put(65,110){ {\large \bf Orbit IX} }
% \put(265,110){ {\large \bf Orbit X} }
% \put(155,-10){ {\large \bf Orbit XI} }
% \end{picture}
% \centerline{\includegraphics[scale=.42]{page2.eps}}
% \label{figure:list2}
% \end{figure}

\newpage
\noindent We conclude by listing some abelian cycles which indicate how to find more sums of index-matched basis elements.
\vskip .5in
\begin{figure}[htpb!]
\begin{picture}(0,0)(0,11)
%\put(65, 340){  {\large \bf LL2} }
%\put(265, 340){  {\large \bf CL1} }
%\put(110, 170){ {\large \bf 3-term sum from QQ7}}
%\put(110, -10){ {\large \bf 2-term sum from QQ8}}
\end{picture}
\centerline{\includegraphics[scale=.9]{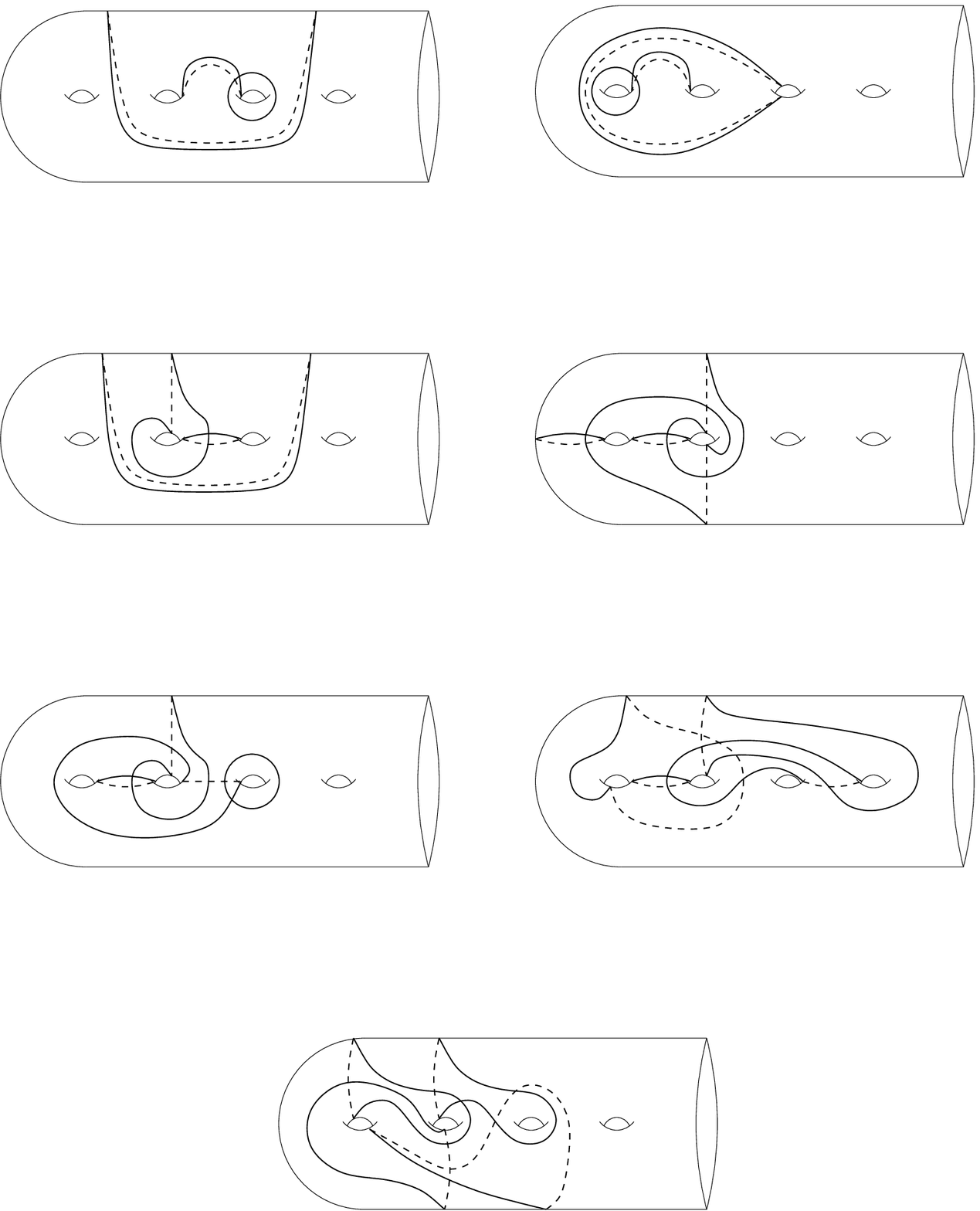}}
\label{figure:indexmatched}
\end{figure}

%\newpage

\noindent
Tara E. Brendle:\\
Dept. of Mathematics, Louisiana State University\\
Lockett Hall\\
Baton Rouge, LA  70803\\
E-mail: brendle@math.lsu.edu
\medskip

\noindent
Benson Farb:\\
Dept. of Mathematics, University of Chicago\\
5734 University Ave.\\
Chicago, Il 60637\\
E-mail: farb@math.uchicago.edu


\begin{thebibliography}{ABCDEF}
\small

\bibitem[Ak]{Ak}
T. Akita, Homological infiniteness of Torelli groups, {\em Topology},
Vol. 40 (2001), no. 2, 213--221.

\bibitem[AM]{AM}S. Akbulut and J.~McCarthy,
{\em Casson's invariant for oriented homology 3-spheres: an exposition},
Mathematical Notes 36, Princeton University Press, 1990.

\bibitem[BC]{BC}
J. Birman and R. Craggs, The $\mu$-invariant of 3-manifolds and certain
structural properties of the group of homeomorphisms of a closed,
oriented 2-manifold, {\em Trans. of the AMS}, Vol. 237 (1978),
p. 283-309.

\bibitem[BF]{BF}
D. Biss and B. Farb, ${\cal K}_g$ is not finitely generated, {\em
Invent. Math.}, Vol. 163, pp. 213--226 (2006).

\bibitem[Br]{Br}
K. Brown, {\em Cohomology of groups}, Graduate Texts in Mathematics,
Vol. 87, Springer, 1982.

\bibitem[Ga]{Ga}
S.~Galatius, Mod $p$ homology of the stable mapping class group, {\em Topology}, Vol. 43 (2004), pp.1105-1132.

\bibitem[GS]{GS}
R. Gompf and A. Stipsicz, 4-Manifolds and Kirby Calculus, AMS
Graduate Studies in Mathematics Vol. 20, 1999.

\bibitem[Hai]{Hai}
R. Hain, Infinitesimal presentations of the Torelli groups, {\em JAMS}, Vol. 10, No. 3 (1997), pp. 597-651.

\bibitem[Ha]{Ha}
J. Harer, The second homology group of the mapping class group of
an orientable surface, {\em Invent. Math.}, Vol. 72 (1983), pp.
221-239.

\bibitem[Hat]{Hat}
A. Hatcher, {\em Algebraic Topology}, Cambridge University Press, 2002.

\bibitem[Jo1]{Jo1}
D. Johnson, Quadratic forms and the Birman-Craggs homomorphisms,
{\em Trans. of the AMS}, Vol. 261, No. 1 (1980), pp. 423-422.

\bibitem[Jo2]{Jo2}
D. Johnson, The structure of the Torelli group I: A finite set of
generators for $\cal{I}$, {\em Annals of Math.}, Vol. 118, No. 3
(1983), pp. 423-422.

\bibitem[Jo3]{Jo3}
D. Johnson, The structure of the Torelli group II: A characterization of the group generated by twists on bounding curves, {\em Topology}, Vol. 24, No. 2 (1985), pp. 113-126.

\bibitem[Jo4]{Jo4}
D. Johnson, The structure of the Torelli group III: the abelianization of $\T$, {\em Topology}, Vol. 24, No. 2 (1985), pp. 127-144.

\bibitem[Jo5]{Jo5}
D. Johnson, A survey of the Torelli group, {\em Contemporary
Mathematics}, Vol. 20 (1983), pp. 165-179.

\bibitem[Jo6]{Jo6}
D. Johnson, An abelian quotient of the mapping class group $\T_g$, {\em Math. Ann.}, Vol. 249 (1980), pp. 225-242.

\bibitem[MW]{MW}
I.~Madsen and M.~Weiss, The stable moduli space of Riemann surfaces:  Mumford's conjecture, preprint (2002), arXiv:math.AT/0212321.

\bibitem[MM]{MM}
D. McCullough and A. Miller, The genus 2 Torelli group is not finitely generated, {\em Topology Appl.}, Vol. 22 (1986), pp. 43-49.

\bibitem[Me]{Me}
G. Mess, The Torelli groups for genus 2 and 3 surfaces, {\em Topology}, Vol. 31, No. 4 (1992), pp. 775-790.

\bibitem[Mo1]{Mo1}
S. Morita, Casson's invariant for homology 3-spheres and
characteristic classes of surface bundles I, {\em Topology}, Vol.
28, No. 3 (1989), pp. 305-323.

\bibitem[Mo2]{Mo2}
S. Morita, On the structure of the Torelli group and the Casson invariant, {\em Topology}, Vol 30, No. 4 (1991), pp. 603-621.

\bibitem[Ra]{Ra} V. Ravikumar, private communication.

\bibitem[Sa]{Sa}
T. Sakasai, The Johnson homomorphism and the third rational cohomology group of the Torelli group, UTMS preprint 2003-21, (http://kyokan.ms.u-tokyo.ac.jp/users/preprint/ps/2003-21.ps).

\bibitem[vdB]{vdB}
B. van den Berg, On the abelianization of the Torelli group,
Ph. D. Thesis, Universiteit Utrecht (2003),
(https://dspace.library.uu.nl:8443/retrieve/25176/inhoud.htm).

\end{thebibliography}
\end{document}